\documentclass[12pt]{amsart}
\usepackage{amssymb}
\usepackage{mathptmx}
\usepackage[all]{xy}

\newcommand{\Z}{\mathbb{Z}}
\newcommand{\R}{\mathbb{R}}

\newcommand{\Q}{\mathbb{Q}}

\newcommand{\Dc}{\mathcal{D}}

\DeclareMathOperator{\chara}{char}

\DeclareMathOperator{\codim}{codim}
\DeclareMathOperator{\depth}{depth}
\DeclareMathOperator{\sdepth}{sdepth}

\DeclareMathOperator{\lk}{lk_{\Delta} }

\DeclareMathOperator{\supp}{supp}

\DeclareMathOperator{\fc}{\mathfrak{c}}

\DeclareMathOperator{\pnt}{\raise 0.5mm \hbox{\large\bf.}}

\newtheorem{thm}{\bf Theorem}[section]

\newtheorem{cor}[thm]{\bf Corollary}

\newtheorem{conj}[thm]{\bf Conjecture}
\newtheorem{quest}[thm]{\bf Question}
\newtheorem{prob}[thm]{\bf Problem}

\newtheorem{lemdef}[thm]{\bf Lemma and Definition}

\theoremstyle{definition}
\newtheorem{defn}[thm]{\bf Definition}

\newtheorem{rem}[thm]{\bf Remark}
\newtheorem{ex}[thm]{\bf Example}

\title{Glicci simplicial complexes}

\author[Uwe Nagel]{Uwe Nagel${}^{*\, \dagger}$}
\address{Department of Mathematics, University of Kentucky, 715 Patterson Office Tower, Lexington, KY 40506-0027, USA}
\email{uwenagel@ms.uky.edu}

\author{Tim R\"omer}
\address{Institut f\"ur Mathematik, Universit\"at Osna\-br\"uck, Albrechtstr. 28a, 49069 Osna\-br\"uck,  Germany}
\email{troemer@uos.de}

\thanks{${}^*$ Part of the work for this paper was done while the first
author was sponsored by the National Security Agency under Grant
Number H98230-07-1-0065.}
\thanks{$^{\dagger}$ Current address: Institute for Mathematics \& its Applications,
University of Minnesota, Minneapolis, MN 55455, USA}

\begin{document}

\begin{abstract}
One of the main open questions in liaison theory is
whether every homogeneous Cohen-Macaulay ideal in a polynomial ring is glicci, i.e.\
if it is in the G-liaison class of a complete intersection.
We give an affirmative answer to this question for Stanley-Reisner ideals
defined by simplicial complexes that are weakly vertex-decomposable.
This class of complexes includes matroid, shifted and Gorenstein complexes respectively.
Moreover, we construct a simplicial complex
which shows that the property of being glicci depends on the characteristic of the base field. 
As an application of our methods we establish new evidence for two conjectures of Stanley on partitionable complexes and Stanley decompositions. 
\end{abstract}


\maketitle


%
%
%
\section{Introduction}

Liaison theory provides an equivalence relation among equidimensional subschemes of fixed dimension. It has found numerous  applications in Algebraic Geometry and Commutative Algebra. Here we adapt its methods to investigate simplicial complexes. Our studies also give rise to interesting arithmetic questions. 

Let us recall some definitions of this theory.
Let $S=K[x_1,\dots,x_n]$ be a polynomial ring over an infinite field $K$.
Two homogeneous ideals $I,J \subset S$ are CI-linked
by a complete intersection $\fc$ if $\fc : I= J$ and $\fc : J= I$.
The transitive closure of this operation leads to the concept of
CI-liaison classes.
CI-liaison theory is well-understood in codimension 2 (see, e.g.,  \cite{migbook}). 
For example, a result of Gaeta \cite{GA50} (see \cite{PESZ74} for a modern proof) implies that every
homogeneous codimension 2 Cohen-Macaulay ideal $I \subset S$ is {\em licci},
i.e.\
it is in the CI-\underline{li}aison \underline{c}lass of a \underline{c}omplete \underline{i}ntersection.

In codimension 3 the situation becomes already much more complicated.
It is known that not every codimension 3 Cohen-Macaulay ideal is licci
(see, e.g.,  Huneke-Ulrich \cite{HU-Annals}).
This is one motivation to link with Gorenstein ideals instead of complete intersection ideals. This idea has been developed to the theory of G-liaison that is coarser than CI-liaison though several important properties generalize (see \cite{KMMNP, migbook, N-gorliaison, schenzel} for details).
One of the main open problems in G-liaison theory is:
\begin{quest}[\cite{KMMNP}]
\label{bigquest}
Is every homogeneous Cohen-Macaulay ideal in $S$ {\em glicci}, i.e., 
is it in the G-\underline{li}aison \underline{c}lass of a \underline{c}omplete \underline{i}ntersection?
\end{quest}
Several classes of Cohen-Macaulay ideals that are of interest in Algebraic Geometry or Commutative Algebra 
are known to be glicci. In this paper we study Question \ref{bigquest} for homogeneous Cohen-Macaulay ideals
which are derived from objects considered in Algebraic Combinatorics.
Recall that
$\Delta$ is called an (abstract) {\em simplicial complex} on
$[n]=\{1,\dots,n\}$ if $\Delta$ is a subset of the power
set of $[n]$ which is
closed under inclusion, i.e.\ if $F \subseteq G$ and $G \in \Delta$,
then $F\in \Delta$.
The elements $F$ of $\Delta$ are called {\em faces}, and the maximal
elements under inclusion are called {\em facets}.

The connection to algebra is provided by the following construction.
For a subset $F
\subset [n]$, we
write $x_F$ for the squarefree monomial $\prod_{i \in F} x_i$.
The {\em Stanley-Reisner ideal} of $\Delta$ is
$I_\Delta =(x_F : F \subseteq [n],\ F\not\in \Delta)$
and the corresponding {\em Stanley-Reisner ring} is $K[\Delta]=S/I_\Delta$.
We will say that $\Delta$
has an algebraic property like Cohen-Macaulayness if 
$K[\Delta]$ has this
property. 
For more details on simplicial complexes, Stanley-Reisner rings and their algebraic properties
we refer to the books of
Bruns-Herzog \cite{BH-book} and Stanley \cite{ST96}.

Given a Cohen-Macaulay complex  $\Delta$, it is natural to ask whether it is glicci, i.e.\
if $I_\Delta$ is a glicci ideal in $S$. Since we are interested in squarefree monomial ideals, we 
study the slightly stronger  property of being {\em squarefree glicci} (see Definition \ref{def-sq-glicci})
which implies being glicci, but is naturally defined in the context of simplicial complexes.
In general we can not answer  Question \ref{bigquest}. But for several 
classes of simplicial complexes we do give an affirmative answer.

In the following assume that $\{i : \{i\} \in \Delta \}=[n]$.
A  recent result of Casanellas-Drozd-Hartshorne \cite{CDH} says 
that each homogeneous Gorenstein ideal $I\subset S$ is glicci.
The proof is complicated and non-constructive.
If $\Delta$ is a Gorenstein complex, then we even show that $\Delta$ is squarefree glicci.
Recall that a complex $\Delta$ is called a {\em matroid} if, for all $W \subseteq [n]$, the restriction $\Delta_W=\{F \in \Delta: F \subseteq W \}$ is a pure simplicial complex. 
We show that in this situation $\Delta$ is squarefree glicci.
Analogously to the idea that a generic initial ideal of a given homogeneous ideal $I$ in $S$
can be used to study algebraic properties of $I$,
one can associate to every simplicial complex a shifted simplicial complex
and these complexes share many combinatorial properties 
(see, e.g., \cite{HE02} or \cite{KA01} for details). 
Recall that a simplicial complex $\Delta$ is called {\em shifted}
if for all $F \in \Delta$, $j \in F$ and $j<i$ such that $i \not\in F$ we have
$F-\{j\} \cup \{i \} \in \Delta$. We prove that every shifted simplicial complex is squarefree glicci.
Hence, for Cohen-Macaulay simplicial complexes, the answer to \ref{bigquest} is affirmative 
up to ``shifting.'' This result is the combinatorial counterpart of one of the main results in \cite{MN4}.
In order to prove that $\Delta$ is squarefree glicci in each of the cases mentioned above, 
we introduce the notion of weakly vertex-decomposable simplicial complexes 
which is a weaker property than the well-known property of being vertex-decomposable.
Our main result establishes that 
weakly vertex-decompo\-sable simplicial complexes are squarefree glicci.
Then it suffices to observe that matroid, shifted and Gorenstein complexes respectively are
weakly vertex-decompo\-sable. 

The concept of weakly vertex-decomposable simplicial complexes is introduced in Section \ref{sec-comb}. There we also prove our main results mentioned above. As prerequisite we  discuss the relevant parts of G-liaison theory in Section \ref{sec-liaison}. Combinatorial applications of our methods are considered in  Section \ref{decompositions}. 

A fundamental open problem for simplicial complexes is a conjecture of Stanley
(see \cite[Problem 6]{ST00})
which states that a Cohen-Macaulay simplicial complex is partitionable.
Using ideas from liaison theory
we prove that if $\Delta$ is
2-CM  and $I_\Delta \subset S$ has codimension 3,
then the  complex $\Delta$ is partitionable.
In particular, this extends the recent result of
Herzog-Jahan-Yassemi \cite{HEJAYA07}
that $\Delta$ is partitionable if
$I_\Delta$ is Gorenstein of codimension 3. Our generalization is in fact obtained by showing that if $\Delta$ is a  2-CM complex  such that $I_{\Delta}$ has codimension 3, then $K[\Delta]$ has a suitable Stanley decomposition. This proves Stanley's conjecture \cite{ST82} about such decompositions in this particular case. 

In the last part of the paper we discuss several 
properties of pure simplicial complexes that are considered in or are related to this paper and their relationships.
In particular, we observe that the triangulation of the real projective plane $\mathbb{P}^2$
as given in \cite[page 236]{BH-book} is not weakly vertex-decomposable.
Thus this complex could be a candidate for which the answer to Question \ref{bigquest} is negative.
Moreover, we construct a simplicial complex which is
weakly vertex-decomposable if $\chara K \neq 2$
and thus glicci, but is not glicci
if $\chara K =2 $. This demonstrates the remarkable fact that the property of being glicci does depend on the characteristic
of the base field. 

%
%
%
\section{Basic double links}
\label{sec-liaison}

We recall results and concepts from Gorenstein liaison theory and specialize them to the case of simplicial complexes.

Two homogeneous ideals $I, J \subset S = K[x_1,\ldots,x_n]$ are said to be {\em G-linked} (in one step) by the Gorenstein ideal $\fc \subset S$ if
\[
\fc : I = J \quad \mbox{and} \quad \fc : J = I.
\]
In this case we write $I \sim_{\fc} J$.
Note that this implies that the ideals $I$ and $J$ are unmixed and have the same codimension as $\fc$. The concept of Gorenstein liaison is obtained if one takes the transitive closure of the above operation, that is, $I$ and $J$ are in the same G-liaison class if and only if there are Gorenstein ideals $\fc_1,\ldots,\fc_s$ such that
\[
I = I_0 \sim_{\fc_1} I_1 \sim_{\fc_2} \cdots \sim_{\fc_s} I_s = J.
\]
If we insist that all the Gorenstein ideals $\fc_1,\ldots,\fc_s$ are in fact complete intersections, then we get the more classical concept of liaison. We will refer to it here as CI-liaison. However, it is crucial for this paper to consider the more general G-liaison. For extensive information on it we refer to \cite{migbook} and \cite{MN-torino}.

Of particular interest are the equivalence classes that contain a complete intersection. We say that the ideal $I$ is {\em glicci} if it is in the {\underline G}-\underline{li}aison {\underline c}lass of a {\underline c}omplete {\underline i}ntersection. It is {\em licci} if it is in the CI-\underline{li}aison {\underline c}lass of a {\underline c}omplete {\underline i}ntersection.
G-liaison is much more flexible than CI-liaison (see, e.g., Example \ref{ex-not-sq-free-glicci}), and this is important for our investigations.
Notice also that every glicci ideal is Cohen-Macaulay and that all complete intersections of the same codimension are in the same CI-liaison class.

We want to use G-liaison to study simplicial complexes by applying it to their Stanley-Reisner ideals. Abusing terminology, we say that simplicial complexes on $[n]$ are linked if their Stanley-Reisner ideals have this property. While linking ideals we would like to stay within the class of squarefree monomial ideals as much as possible. However, in general even for arbitrary monomial ideals it is too restrictive to require that all the ideals involved in the links are monomial.

\begin{ex}
\label{ex-not-sq-free-glicci}
Consider the ideal $I = (x_1,x_2,x_3)^2 \subset K[x_1,x_2,x_3]$.  It  has a linear free resolution.  By \cite{HU-Annals} this implies that $I$ is not licci.
However, $I$ is glicci (see Example \ref{ex-glicci}(ii)), but it is not possible to G-link $I$ to a complete intersection using only monomial ideals. Indeed, this follows because each artinian Gorenstein monomial ideal of codimension 3 is a complete intersection. For example, this can be deduced from the structure theorem for monomial Gorenstein ideals of codimension 3  established in \cite{BH-95}.
\end{ex}

In view of this example we require only that every other ideal is monomial.

\begin{defn}
\label{def-sq-glicci} \ 

(i) The squarefree monomial ideal $I \subset S$ is said to be {\em squarefree glicci} if  there is a chain of links in $S$
\[
I = I_0 \sim_{\fc_1} I_1 \sim_{\fc_2} \cdots \sim_{\fc_{2s}} I_{2s},
\]
where $I_j$ is a squarefree monomial ideal whenever $j$ is even and $I_{2s}$ is a complete intersection.

(ii) Let $\Delta$ be simplicial complex  with existing vertices $\{i : \{i\} \in \Delta\} = [n]$. Then $\Delta$ is called {\em squarefree glicci} if  $I_{\Delta} \subset S$ has this property.
\end{defn}

Let $I \subset S$ be an ideal and let $R = S[y]$ be the polynomial ring over $S$ in the variable $y$. If $I$ is  glicci, then so is the extension ideal $I \cdot R$ because the links in $S$ also provide links in $R$. This implies in particular, that if $\Delta$ is squarefree glicci, then so is any cone over $\Delta$. We will use this fact frequently in this note.


There is a simple construction that often allows us to link a given simplicial complex in two steps to a subcomplex.

\begin{lemdef}
\label{bdl}
Let $\fc \subset J \subset S$ be squarefree monomial ideals and let $x_k \in S$ be a variable that does not divide any minimal monomial neither in $J$ nor in $\fc$. If $\fc$ is Cohen-Macaulay and $J$ is unmixed such that $\codim J = \codim \fc + 1$, then $I := x_k J + \fc$ is a squarefree monomial ideal that is G-linked in two steps to $J$. We say that $I$ is a  {\em basic double link} of $J$ on $\fc$.
\end{lemdef}

\begin{proof}
This follows from \cite[Proposition 5.10]{KMMNP}. For the convenience of the reader we sketch the argument in this special case. Write $\bar{\mbox{\ }}$ for the images in $\bar{S} = S/\fc$. Since every squarefree monomial ideal is locally a complete intersection, the canonical module of $\bar{S}$ is (up to a degree shift) isomorphic to an ideal $\omega$ of height one in $\bar{S}$ (see, e.g., \cite[Proposition 3.3.18]{BH-book}). This means that there is a homogeneous Gorenstein ideal $G' \subset S$ containing $\fc$ such that $G'/\fc = \omega$. Since $\codim J > \codim \fc$, there is a homogeneous polynomial $f \in J$ such that $\bar{f}$ is $\bar{S}$-regular. But $\bar{f} \omega$ is also (up to a degree shift) isomorphic to the canonical module of $\bar{S}$, hence $G :=  f \cdot G' + \fc$ is a Gorenstein ideal contained in $J$. Similarly, we get that $ x_k G + \fc = x_k f \cdot G' + \fc$ is a Gorenstein ideal contained in $I$ because $x_k$ is also $\bar{S}$-regular. The latter fact also implies  that $(\bar{x_k} \bar{f} \cdot \omega) : (\bar{x_k} \cdot \bar{J}) = (\bar{f} \cdot \omega) : \bar{J}$ in $\bar{S}$. Back in $S$ this means $I' := (x_k G + \fc ) : (x_k J + \fc) = G : J$, which provides the G-links $I \sim_{(x_k G + \fc)} I'
 \sim_{G} J$, as claimed.
\end{proof}

\begin{rem}
\label{rem-bdl}
\ 

(i) The above definition is a very special case of the more general concept of a basic double link for G-liaison as introduced in \cite{KMMNP}. For previous uses of basic double links we refer to, for example, \cite{HU}, \cite{KMMNP},  \cite{MN4},  \cite{MN3}, \cite{MN-tetra}. 

(ii) If $I$ is a basic double link of $J$, then $I$ is Cohen-Macaulay if and only if $J$ has this property (see, e.g., \cite{migbook}).

(iii) Let $\Delta$ be a simplicial complex on $[n]$. Recall that each  $F \subseteq [n]$ induces the  following simplicial subcomplexes of $\Delta$:
the {\em link of $F$}
$$
\lk F = \{G \in \Delta : F \cup G \in \Delta, F \cap G = \emptyset
\},
$$
and the {\em deletion}
$$
\Delta_{-F} = \{G \in \Delta : F \cap G = \emptyset  \}.
$$

Consider any $k \in [n]$. 
Then the cone over the link $\lk k$ with apex $k$ considered as complex on $[n]$ has as Stanley-Reisner ideal $J_{\lk k} = I_{\Delta} : x_k$,  and the Stanley-Reisner ideal of the deletion $\Delta_{-k}$ considered as a complex on $[n]$ is $(x_k, J_{\Delta_{-k}})$  where $J_{\Delta_{-k}} \subset S$ is the extension ideal of the  Stanley-Reisner ideal of $\Delta_{-k}$ considered as a  complex on $[n] \setminus \{k\}$. Note that $x_k$ does not divide any of the minimal generators of $J_{\Delta_{-k}}$, thus $x_k$ is not a zerodivisor on $S/J_{\Delta_{-k}}$. Furthermore, it follows that $I_{\Delta} = x_k J_{\lk k} +  J_{\Delta_{-k}}$. Hence, if $\Delta$ is pure and if the deletion $\Delta_{-k}$ is Cohen-Macaulay and has the same dimension as $\Delta$, then $\Delta$ is a basic double link of the cone over its link $\lk k$, where both are considered as complexes on $[n]$.  Moreover, each such basic double link provides the exact sequence
\begin{equation*} \label{eq-seqqq}
0
\to
S/J_{\lk k} (- \deg x_k)
\overset{x_k}{\to}
S/I_\Delta
\to
S/I_{\Delta_{-k}}
\to
0.
\end{equation*}
\end{rem}

\begin{ex}
\label{ex-glicci} \  

(i) Consider the simplicial complex on $[4]$ consisting of 4 vertices. Its Stanley-Reisner ideal is
 $I_{\Delta}  = (x_1 x_2, x_1 x_3, x_1 x_4, x_2 x_3, x_2 x_4, x_3 x_4)$. It  has a linear free resolution, from which it follows by \cite{HU-Annals} that it is not licci.
However, $I$ is squarefree glicci because
\[
I = x_4 \cdot (x_1, x_2, x_3) + (x_1 x_2, x_1 x_3, x_2 x_3)
\]
provides that $I$ is a basic double link of $(x_1, x_2, x_3)$.

(ii) Similarly, using the general version  of basic double linkage (\cite[Proposition 5.10]{KMMNP}), 
$(x_1,x_2,x_3)^2 = x_3 \cdot (x_1,x_2,x_3) + (x_1,x_2)^2$ shows that $(x_1,x_2,x_3)^2$ is glicci.
\end{ex}

%
%
%

\section{Weakly vertex-decomposable complexes} \label{sec-comb}

The goal of this section is to identify a combinatorially defined class of simplicial complexes that consists of  squarefree glicci complexes. It includes, for example, shifted, matroid,  and Gorenstein complexes.

Following \cite{PB-80} (see also \cite[Definition 11.1]{BW97}), a pure simplicial complex $\Delta$ is said to be {\em vertex-decomposable} if $\Delta$ is a simplex or equal to $\{\emptyset\}$, or there exists a vertex $k$ such that $\lk k$ and $\Delta_{-k}$ are both pure and vertex-decomposable and $\dim \Delta = \dim \Delta_{-k} = \dim \lk k + 1$.

We now propose a less restrictive concept that is also defined  recursively:

\begin{defn} \label{def-CM-vd} The pure simplicial simplex $\Delta \neq \emptyset$ on $[n]$ is said
  to be {\em weakly vertex-decompo\-sable} if
  there is some $k \in [n]$ such that $\Delta$ is a cone over the
  weakly vertex-decomposable deletion $\Delta_{-k}$ or there is some $k \in [n]$ such that $\lk k$ is weakly vertex-decomposable and
  $\Delta_{-k}$ is Cohen-Macaulay of the same dimension as
  $\Delta$.
\end{defn}

Observe, that if $\Delta$ is not a cone over $\Delta_{-k}$,
then $\Delta_{-k}$ and $\Delta$ have automatically the same
dimension.

\begin{ex}
\label{ex-CM-vd}
Simplicial complexes on $[n]$ of high dimension are often weakly vertex-decompo\-sable:

(i) If $\dim \Delta = n-1$, then $\Delta$ is a simplex, thus it is  weakly vertex-decompo\-sable.

(ii) If $\Delta$ is pure of dimension $n-2$, then for any vertex $\{k\} \in \Delta$, the Stanley-Reisner ideals of $\Delta$ and the cones over $\lk k$ and $\Delta_{-k}$, respectively,  are principal ideals, so $\Delta$ is weakly vertex-decomposable and Cohen-Macaulay.
\end{ex}

\begin{thm}
\label{thm-CM-vd}
If $\Delta$ is weakly vertex-decomposable, then $\Delta$ is squarefree glicci. In particular, $\Delta$ is Cohen-Macaulay.

\end{thm}

\begin{proof}
We use induction on  $n$ to show that there is a finite sequence of basic double links starting with $\Delta$ that ends with a complete intersection.

If  $n = \dim \Delta + 1$, then the Stanley-Reisner ring of $\Delta$ is regular, thus in particular a complete intersection.
Let $n \geq \dim \Delta + 2$. If there is a $k \in [n]$ such that $\Delta$ is a cone over the weakly vertex-decomposable deletion $\Delta_{-k}$, then $\Delta_{-k}$ considered as complex over $[n] \setminus \{k\}$ is squarefree glicci by induction. As discussed above, this implies that $\Delta$ is squarefree glicci.
 
It remains to consider the case, where there is a vertex $k$ of $\Delta$ such that $\lk k$ is weakly vertex-decomposable and $\Delta_{-k}$ is Cohen-Macaulay of the same dimension as $\Delta$. Then we know by Remark \ref{rem-bdl}(iii) that $\Delta$ is a basic double link of the cone over $\lk k$. By induction, $\lk k$ considered as complex on $[n] \setminus \{k\}$ is squarefree glicci, thus so is the cone over it. It follows that $\Delta$ is squarefree glicci.
\end{proof}

In \cite{MN4} it has been shown that each Cohen-Macaulay strongly stable ideal is glicci. Our first consequence is that the squarefree analogue is true as well. 

\begin{cor}
\label{cor-shifted-gl}
Each Cohen-Macaulay shifted complex  is squarefree glicci.
\end{cor}

\begin{proof}
According to \cite{BW85} each Cohen-Macaulay shifted complex is vertex-decomposable, thus we conclude by Theorem \ref{thm-CM-vd}.
\end{proof}

Notice that the passage from a simplicial complex $\Delta$ to its symmetric shift $\Delta^s$  is analogous to the passage from a homogeneous ideal to its generic initial ideal (see \cite{ARHEHI98}).  Furthermore, if $\Delta$ is Cohen-Macaulay, then so is $\Delta^s$. In this sense, Corollary \ref{cor-shifted-gl} shows that every Cohen-Macaulay complex is squarefree glicci up to ``shifting.''
\smallskip

The following result concerns matroids.

\begin{cor}
Each matroid is squarefree glicci.
\end{cor}

\begin{proof}
Let $k \in [n]$. Then it is well-known that $\lk k$ and $\Delta_{-k}$ are the corresponding link and deletion in the sense of matroid theory. In particular, they are again matroids (see, e.g.,  \cite{Oxley-book-92}). Hence it follows by induction on the number of vertices that each matroid is vertex-decomposable. This completes the argument.
\end{proof}

Following \cite{Bacl-82}, the complex $\Delta$ is said to be {\em $2$-CM} or doubly Cohen-Macaulay if, for each existing vertex $\{k\} \in \Delta$, the deletion $\Delta_{-k}$ is Cohen-Macaulay of the same dimension as $\Delta$.

\begin{cor}
\label{cor-2cm}
Each 2-CM complex is squarefree glicci.
\end{cor}

\begin{proof}
In order to see that each $2$-CM complex is weakly vertex-decompo\-sable it suffices to check that its link with respect to any vertex is again $2$-CM. This has been shown in \cite{Bacl-82}, \cite{Bacl-83}.
\end{proof}

A recent result by Casanellas-Drozd-Hartshorne (see \cite{CDH}) says that each Gorenstein ideal is glicci. The proof is non-constructive and relies on the theory developed in \cite{CDH}. Thus, it is somewhat surprising that our method provides an even stronger result for Stanley-Reisner ideals.

\begin{cor}
\label{cor-Gor}
Each simplicial homology sphere is squarefree glicci.
\end{cor}

\begin{proof}
Note that the Stanley-Reisner ring of a homology sphere is Gorenstein. Furthermore,
Hochster's Tor formula provides that each Gorenstein ideal is $2$-CM (see \cite{Bacl-82}).
\end{proof}

\begin{rem}
Actually, the proofs of the above results establish a stronger result that includes some monotonicity.  Indeed, we show in all the cases above that if the complex is not a simplex, then it is a basic double link of a proper subcomplex. It terms of liaison theory, this means that the schemes defined by our squarefree glicci ideals can be G-linked to a complete intersection by descending basic double links.
\end{rem}

We conclude this section with an observation about the relevance of the characteristic of the ground field.

\begin{rem}
Notice that the property of being
squarefree glicci and even being glicci depends on the characteristic of the ground field (see Example \ref{ex-char-dep}).

However, the results in \cite{HU} imply that for monomial ideals $I \subset S$ that are artinian, i.e.\ $S/I$ is a finite-dimensional $K$-vector space, the property of being licci is independent of the characteristic of the ground field.
\end{rem}

This suggests the following 

\begin{prob}
\label{prob-char}
Is for squarefree monomial ideals the property of being licci independent of the characteristic of the ground field?

Furthermore, motivated by Example \ref{ex-char-dep}, one might also ask whether there is a squarefree monomial ideal $I \subset K[x_1,\ldots,x_n]$ and a prime number $p$ such that $I$ is licci if the characteristic of $K$ is not $p$, but $I$ is not licci if $\chara K = p$?
\end{prob}

%
%

\section{Stanley Decompositions}
\label{decompositions}

Let $R$ be any standard graded $K$-algebra over an infinite field $K$, i.e., 
$R$ is a finitely generated graded algebra $R=\bigoplus_{i \geq 0} R_i$
such that $R_0=K$ and $R$ is generated by $R_1$.
There are several characterizations of the depth of such an algebra.
We use the one that $\depth R$ is the
maximal length of a regular $R$-sequence consisting of linear forms.

The goal of this section is to make a contribution to a conjecture of  Stanley in \cite{ST82}.
We follow the terminology of Herzog-Jahan-Yassemi \cite{HEJAYA07} who proved this conjecture
in several cases.
For this let $x_F=\prod_{i \in F}x_i$ be a squarefree monomial for some $F \subseteq [n]$
and $Z \subseteq \{x_1,\dots,x_n\}$.
The $K$-subspace $x_F K[Z]$ of $S = K[x_1,\dots,x_n]$ is the subspace generated by monomials $x_F u$,
where  $u$ is a monomial in the polynomial ring $K[Z]$. It 
is called a {\em squarefree Stanley space} if $\{x_i : i \in F\} \subseteq Z$.
The {\em dimension} of this Stanley space is $|Z|$.
Let $\Delta$ be a simplicial complex on $[n]$.
A {\em squarefree Stanley decomposition}
$\Dc$ of $K[\Delta]$ is a finite direct sum $\bigoplus_{i} u_i K[Z_i]$ of squarefree
Stanley spaces which is isomorphic as a $\Z^n$-graded $K$-vector space to $K[\Delta]$, i.e.\ 
$$
K[\Delta] \cong \bigoplus_{i} u_i K[Z_i].
$$
We denote by $\sdepth \Dc$  the minimal dimension of a Stanley space in $\Dc$ and we define
$$
\sdepth K[\Delta] =\max \{ \sdepth \Dc : \Dc \text{ is a Stanley decomposition of } K[\Delta]\}.
$$
Stanley conjectured in \cite[Conjecture 5.1]{ST82}
 the following upper bound for the depth of $K[\Delta]$.
\begin{conj}[Stanley]
\label{conjdec}
If $\Delta$ is a simplicial complex on  $[n]$, then
\[
\depth K[\Delta] \leq \sdepth K[\Delta].
\]
\end{conj}
Observe that in Stanley's original conjecture a more general situation is considered.
However, in \cite[Theorem 3.3]{HEJAYA07} it was shown that the conjecture stated as above
is equivalent to Stanley's conjecture in the case of Stanley-Reisner rings.
In this section we prove Conjecture \ref{conjdec} in a special case
using ideas from  liaison theory considered in previous sections of this paper.

For this we recall some notation.
Let $\Delta$ be a simplicial complex on $[n]$ and let $k \in [n]$.
In the polynomial ring $S$
we consider the ideals
$$
J_{\lk k} =( x_F : F \subseteq [n] \setminus \{k\}, F \not\in \lk k )
\text{ and }
J_{\Delta_{-k}} =( x_F : F \subseteq [n] \setminus \{k\}, F \not\in \Delta_{-k}).
$$
If we let
$\lk k$
and
$\Delta_{-k}$
be simplicial complexes on the vertex set $[n]$,
then it follows from the definitions that the corresponding Stanley-Reisner ideals of these complexes on $[n]$ are exactly
$$
I_{\lk k}
=
(x_k)
+
J_{\lk k}
\subset S
\text{ and }
I_{\Delta_{-k}}
=
(x_k)
+
J_{\Delta_{-k}}
\subset S.
$$
Notice that the ideals
$J_{\lk k}$ and $J_{\Delta_{-k}}$
are the Stanley-Reisner ideals of
the cones  with apex $k$ over the simplicial complexes
$\lk k$ and $\Delta_{-k}$ considered as complexes over $[n]\setminus \{k\}$.

Herzog-Jahan-Yassemi \cite{HEJAYA07} proved Conjecture \ref{conjdec}
in the cases that $I_\Delta$ is Cohen-Ma\-cau\-lay of codimension 2 and
Gorenstein of codimension 3.
We now generalize the latter result.

\begin{thm}
Let $\Delta$ be a 2-CM simplicial complex on $[n]$
such that $I_\Delta \subset S$ has codimension 3.
Then
$$
\depth K[\Delta] \leq \sdepth K[\Delta].
$$
In particular, if
$I_\Delta$ is
Gorenstein of codimension 3, then this inequality is true.
\end{thm}
\begin{proof}
Let $V(\Delta)=\{ i \in [n] : \{i\} \in \Delta  \}$ be the set of existing vertices of $\Delta$.
We prove the theorem by induction on $|V(\Delta)|$.

For $|V(\Delta)|=0$ we have
$
S/I_\Delta=S/(x_1,\dots,x_n)=K.
$
Since $\depth K =0$ and $K=1\cdot K$ is a Stanley decomposition, it follows that the inequality $\depth S/I_\Delta \leq \sdepth S/I_\Delta$ is true in this case.
Similarly for $|V(\Delta)|=1$ we have after possibly renumbering the variables that
$
S/I_\Delta=S/(x_2,\dots,x_n)\cong K[x_1].
$
Since $\depth S/I_\Delta=1$ and there exists the Stanley decomposition $1\cdot K[x_1]$
we are also done in this case.

Assume now that $|V(\Delta)|>1$. Let $k \in V(\Delta)$ be an existing vertex of $\Delta$.
The simplicial complex $\lk k$ is a complex on
$[n]$ with $|V(\lk k)|<|V(\Delta)|$ and $ \dim \lk k = \dim \Delta -1$.
It is again 2-CM (see, e.g., \cite{Bacl-82, Bacl-83}).
Let $S'=S/(x_k)$ and $I'_{\lk k}$ be the Stanley-Reisner ideal of
$\lk k$ considered as an complex on $[n]\setminus \{k\}$.
Then $I'_{\lk k}$ is an Cohen-Macaulay ideal of codimension 3 in $S'$
and $S/I_{\lk k} \cong S'/I'_{\lk k}$ as $\Z^n$-graded $K$-algebras.
By induction, the ring $S'/I'_{\lk k}$ has 
a squarefree Stanley decomposition, so  we get an 
isomorphism of $\Z^n$-graded $K$-vector spaces
$$
S'/I'_{\lk k} \cong \bigoplus_{i} u_i K[Z_i']
$$
such that, for all $i$,  $|Z'_i| \geq \depth S'/I'_{\lk k}$ , $x_k \not\in Z_i'$, and $ x_k \nmid u_i$.
Recall that
$I_{\lk k}=(x_k)+J_{\lk k}$ and that $x_k$ does not appear as a factor of
any minimal generator of $J_{\lk k}$.
It follows that
$S/J_{\lk k}$ has the squarefree Stanley decomposition
\begin{equation}
\label{dec_one}
S/J_{\lk k} \cong (S'/I'_{\lk k})[x_k] \cong  \bigoplus_{i} u_i K[Z_i], 
\end{equation}
where $Z_i=Z_i' \cup \{x_k\}$, 
Observe that, \text{ for all } $i$,
$$
|Z_i|
=
|Z'_i|+1
\geq
\depth S/I_{\lk k} +1
=
\depth S/J_{\lk k}.
$$
Next note that the simplicial complex $\Delta_{-k}$ is also a complex on
$[n]$ with $|V(\Delta_{-k})|<|V(\Delta)|$.
We consider again $S'=S/(x_k)$ and let $I'_{\Delta_{-k}}$ be the Stanley-Reisner ideal of
$\Delta_{-k}$ considered as an complex on $[n]\setminus \{k\}$.
Then $I'_{\Delta_{-k}}$ is an Cohen-Macaulay ideal of codimension 2
and  $S/I_{\Delta_{-k}} \cong S'/I'_{\Delta_{-k}}$ as $\Z^n$-graded $K$-algebras.
It follows from \cite[Proposition 1.4]{HEJAYA07}
that the ring $S'/I'_{\Delta_{-k}}$ has
a squarefree Stanley decomposition
and we obtain isomorphisms of $\Z^n$-graded $K$-vector spaces
\begin{equation}
\label{dec_two}
S/I_{\Delta_{-k}} \cong S'/I'_{\Delta_{-k}} \cong \bigoplus_{j} v_j K[Y_i]
\end{equation}
such that $|Y_j| \geq \depth S/I_{\Delta_{-k}}$ for all $j$.
Now we consider the short exact sequence
$$
0
\to
S/J_{\lk k} (-\epsilon_k)
\overset{x_k}{\to}
S/I_\Delta
\to
S/I_{\Delta_{-k}}
\to
0
$$
where $\epsilon_k$ denotes the $k$-th standard basis vector of $\Z^n$.
This sequence together with (\ref{dec_one}) and (\ref{dec_two})  yields
the following decomposition of $S/I_\Delta$ as $\Z^n$-graded
$K$-vector spaces
$$
S/I_\Delta
\cong
\bigoplus_{i} x_k\cdot u_i K[Z_i]
\oplus
\bigoplus_{j} v_j K[Y_i].
$$
Observe that $x_k\nmid u_i$ and $x_k \in Z_i$
and therefore $\supp(x_k\cdot u_i) \subseteq Z_i$ for all $i$.
Hence this decomposition is
a squarefree Stanley decomposition
of $S/I_\Delta$.
Using 
$$
\depth S/I_\Delta
=
\depth S/J_{\lk k}
=
\depth S/I_{\Delta_{-k}}
$$
we get that, \text{ for all } $i,j$,
$$
\depth S/I_\Delta
 \leq |Z_i|
 \text{ and }
\depth S/I_\Delta
 \leq |Y_j|.
$$
Thus
$$
\depth S/I_\Delta \leq \sdepth S/I_\Delta,
$$
and this concludes the proof.
\end{proof}

Next we recall another conjecture of Stanley.
Let $\Delta$ be again a simplicial complex on $[n]$
with facets $G_1,\dots,G_t$. The complex
$\Delta$ is called {\em partitionable} if
there exists a partition $\Delta= \bigcup_{i=1}^t [F_i,G_i]$
where $F_i \subseteq G_i$ are suitable faces of $\Delta$.
Here the interval $[F_i,G_i]$ is the set of faces $\{H \in \Delta : F_i \subseteq H \subseteq G_i\}$.
In \cite[Conjecture 2.7]{ST96} and \cite[Problem 6]{ST00} respectively
Stanley conjectured the following:

\begin{conj}[Stanley]
\label{conjpart}
Let $\Delta$ be a simplicial complex on $[n]$ which is Cohen-Macaulay.
Then $\Delta$ is partitionable.
\end{conj}

Conjecture \ref{conjpart}
is a special case of Conjecture \ref{conjdec}.
Indeed, Herzog-Jahan-Yassemi \cite[Corollary 3.5]{HEJAYA07}
proved that for Cohen-Macaulay simplicial complex $\Delta$ on
$[n]$ we have that $\depth S/I_\Delta \leq \sdepth S/I_\Delta$ if and only if
$\Delta$ is partitionable. As a consequence of this result we obtain:

\begin{cor}
\label{nicepart}
Let
$\Delta$ be a 2-CM simplicial complex on $[n]$
such that  $I_\Delta \subset S$ has codimension 3.
Then Conjecture \ref{conjpart} is true.
\end{cor}

%
%

\section{Examples}

The following diagram displays some properties of pure simplicial complexes and their relationships. 
This list of properties is not complete at all.
We discuss only the ones that are considered in or that are related to this paper.
$$
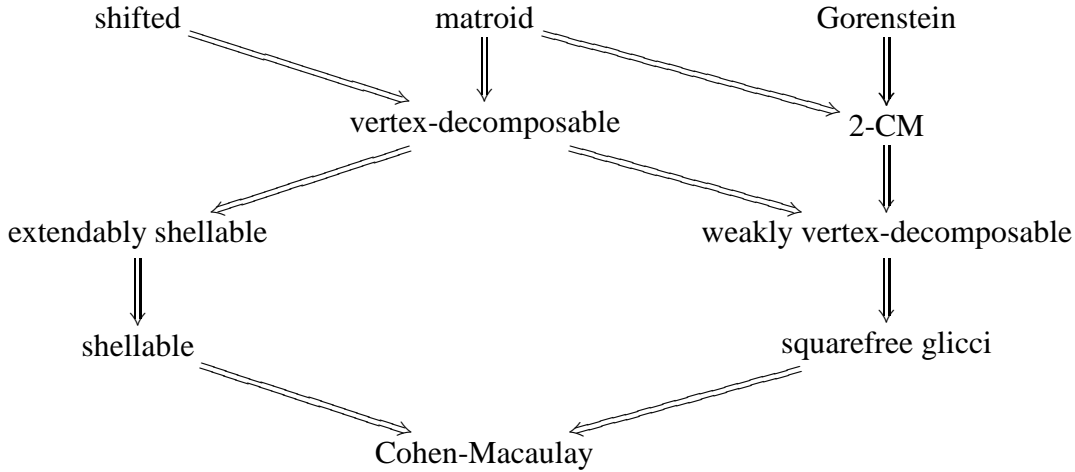
\xymatrix{
\text{shifted}\ar@{=>}[dr]&\text{matroid}              \ar@{=>}[d]\ar@{=>}[dr]                & \text{Gorenstein}\ar@{=>}[d]      \ar@{=>}[d]\\
&\text{vertex-decomposable}  \ar@{=>}[ld]\ar@{=>}[rd]   & \text{2-CM}                       \ar@{=>}[d]\\
\text{extendably shellable}  \ar@{=>}[d]&  & \text{weakly vertex-decomposable} \ar@{=>}[d]\\
\text{shellable}             \ar@{=>}[dr]&               & \text{squarefree glicci}         \ar@{=>}[dl]\\
&\text{Cohen-Macaulay}                  & \\
}
$$
In this section we
discuss (counter)-examples which show that most of the above implications can not be reversed.

We mentioned above that matroids and pure shifted complexes are vertex-decompo\-sable, 
and it is well-known that these implications are strict.
Notice that matroids are also 2-CM.
The converse is not true (see below).

We saw in the proof of \ref{cor-Gor} that a Gorenstein complex is 2-CM.
The simplicial complex in Example \ref{ex-glicci} 
is 2-CM (indeed it is a matroid) and is not Gorenstein. Thus this implication is strict.

Let $\Delta$ be a pure simplicial complex on $[n]$.
The complex $\Delta$ is called {\em shellable}, if
the facets of $\Delta$ can be given a linear order $F_1,\dots,F_t$.
such that $\langle F_i \rangle \cap \langle F_1,\dots, F_{i-1}\rangle$
is generated by a non-empty set of maximal proper faces of $\langle F_i \rangle$
for $2\leq i \leq t$. We call $\Delta$ {\em extendably shellable} if any order of the facets
of $\Delta$ has this property.
Lickorish \cite{L-91} gives examples of simplicial $n$-spheres for $n \geq 3$
that are not shellable (see also \cite{LU04}).
Hence Gorenstein complexes are in general not 
shellable. Thus also the properties 2-CM, weakly vertex-decomposable
and squarefree glicci do not imply shellability or any property above this one.

In the proof of \ref{cor-2cm} we showed that 2-CM complexes are weakly vertex-decomposable.
The converse is not true as the following example shows.

\begin{ex}
Let $S=K[x_1,\dots,x_5]$ and consider the simplicial complex $\Delta$ defined by
$I_\Delta=(x_1x_2, x_1x_3, x_1x_4, x_2x_3, x_2x_4, x_3 x_4 x_5)$.
Then $\Delta$ is shifted of dimension $1$ and $K[\Delta]$
has the Betti-diagram
(e.g.\
use the Eliahou-Kervaire type formula for shifted simplicial complexes
given in \cite{ARHEHI98})
$$
\begin{array}{c|cccc}
  & 0 & 1 & 2 & 3 \\
\hline
0 & 1 & - & - & - \\
1 & - & 5 & 6 & 2 \\
2 & - & 1 & 2 & 1 \\
\end{array}
$$
Thus we see that $K[\Delta]$ is Cohen-Macaulay.
But $K[\Delta]$ is not level (i.e.\ the canonical module is not generated in one degree).
Since 2-CM complexes have Stanley-Reisner rings that are level (see \cite[Theorem 5.7.6]{BH-book}),
the complex $\Delta$ can not be 2-CM. But  each shifted complex is weakly vertex-decomposable.
Therefore
vertex-decomposable and weakly vertex-decomposable complexes, respectively,  are in general not  2-CM.
\end{ex}

The left part of the diagram above
$$
\text{shifted} \Rightarrow
\text{vertex-decomposable} \Rightarrow
\text{extendably shellable} \Rightarrow
\text{shellable}
$$
is well-known and all implications are strict.
See \cite{BW85} for the first implication.
Moriyama-Takeuchi \cite[Theorem B]{MF-03} constructed a simplicial complex
which is extendably shellable, but not vertex-decomposable.
One can construct also examples of shellable, but not extendably shellable  simplicial complexes
(see, e.g.,  \cite[Theorem A]{MF-03}).

A shellable simplicial complex is Cohen-Macaulay.
The standard example of a non-shellable Cohen-Macaulay complex is a triangulation of the real projective plane $\mathbb{P}^2$
on the vertex set $[6]$ (see, e.g., \cite[page 236]{BH-book}).
This triangulation is also not weakly vertex-decomposable as
the following example shows.

\begin{ex}
Let $S=K[x_1,\dots,x_6]$.
Using the notation from \cite[p. 236]{BH-book}, the Stanley-Reisner ideal
of the 
triangulation of the real projective plane $\mathbb{P}^2$
is given by
$$
I_\Delta=
(
x_1x_2x_3,
x_1x_2x_4,
x_1x_3x_5,
x_1x_4x_6,
x_1x_5x_6,
x_2x_3x_6,
x_2x_4x_5,
x_2x_5x_6,
x_3x_4x_5,
x_3x_4x_6)
.
$$
If $\chara K\neq 2$ this is a $2$-dimensional Cohen-Macaulay
complex, while for  $\chara K =2$ this complex is not Cohen-Macaulay
(and thus not shellable).
Assume now that $K=\Q$.
We used Macaulay 2 \cite{MC2} to check that $\Delta$ is not weakly vertex-decomposable.
Indeed, for all $k \in [6]$ the rings $S/J_{\Delta_{-k}}$
are 4 dimensional, but have depth 3 and therefore are not Cohen-Macaulay.
For example, for $k=1$ the corresponding ideal is
$$
J_{\Delta_{-1}}=(
x_2x_3x_6,
x_2x_4x_5,
x_2x_5x_6,
x_3x_4x_5,
x_3x_4x_6).
$$
\end{ex} 
\smallskip 

One of the main open questions in liaison theory is whether every Cohen-Macaulay ideal is glicci.
In view of the above dependence of the Cohen-Macaulayness on the characteristic,  we propose the following: 

\begin{prob} 
Decide whether the Stanley-Reisner ideal of the above triangulation of $\mathbb{P}^2_{\R}$   is glicci. 
\end{prob}

In the next example we show that there exists a extendably shellable simplicial complex
which is not weakly vertex-decomposable.

\begin{ex}
Let $S=K[x_1,\dots,x_6]$.
The ideal of the simplicial complex of \cite[Theorem B, V6F10-6]{MF-03}
is given by
$$
I_\Delta=(x_1x_2x_6,
x_1x_3x_5,
x_1x_4x_5,
x_1x_4x_6,
x_1x_5x_6,
x_2x_3x_4,
x_2x_3x_5,
x_2x_3x_6,
x_2x_4x_6,
x_3x_4x_5)
$$
Moriyama-Takeuchi observed that $\Delta$ is extendable shellable,
but not vertex-decomposable.
We used Macaulay 2 \cite{MC2}
to check that $\Delta$ is also not weakly vertex-decomposable.
\end{ex}

The last example of this section shows that
the two properties being weakly vertex-decomposable and being 
 squarefree glicci  depend on the characteristic of $K$.

\begin{ex}
\label{ex-char-dep}
Let $S = K[x_1,\dots,x_7]$ and consider the ideals
\begin{eqnarray*}
\fc
&=&
(
x_1x_2x_3,
x_1x_2x_4,
x_1x_3x_5,
x_1x_4x_6,
x_1x_5x_6,
x_2x_3x_6,
x_2x_4x_5,
x_2x_5x_6,
x_3x_4x_5,
x_3x_4x_6),\\
J
&=&
(x_1,\dots,x_4),  \text{ and }\\
I&=&x_7 J + \fc.
\end{eqnarray*} 
Notice that $\fc$ is the extension of the Stanley-Reisner ideal of the 
triangulation of the real projective plane $\mathbb{P}^2$ in $6$ variables! Hence, $S/\fc$ is Cohen-Macaulay if and only if $\chara K \neq 2$. Therefore $I$ is a basic double link of the complete intersection $J$ if $\chara K \neq 2$. It follows that in this case $I$ is squarefree glicci and that the induced simplicial complex $\Delta$ is weakly vertex-decomposable. Furthermore, using Macaulay 2
with $K=\Z/31013\Z$ we get the Betti-diagram
$$
\begin{array}{c|ccccc}
  & 0 & 1  & 2  & 3 & 4\\
\hline
0 & 1 & -  & -  & - & -\\
1 & - & 4  & 6  & 4 & 1\\
2 & - & 10 & 25 & 21& 6\\
\end{array}
$$

Now assume that the characteristic of $K$ is 2. Using the exact sequence 
\[
0 \to \fc (-\deg x_7) \to \fc \oplus J (-\deg x_7) \to I \to 0, 
\]
it is not too difficult to check that $S/I$ has depth $2 < \dim S/I = 3$, thus $S/I$ is not Cohen-Macaulay. It follows that $\Delta$ is neither (squarefree) glicci nor weakly vertex decomposable if $\chara K = 2$. 
It is amusing to compare the Betti numbers. Another computation with Macaulay 2
with $K=\Z/2\Z$ gives the Betti-diagram
$$
\begin{array}{c|cccccc}
  & 0 & 1  & 2  & 3 & 4 & 5\\
\hline
0 & 1 & -  & -  & - & - & -\\
1 & - & 4  & 6  & 4 & 1 & -\\
2 & - & 10 & 25 & 21& 7 & 1\\
3 & - & - & -   & 1 & 1 & -\\
\end{array}
$$
\end{ex}

%
%

\end{document}